\documentclass[11pt]{article}
\usepackage[american]{babel}
\usepackage{amsmath}
\usepackage{latexsym}
\usepackage{amssymb}
\usepackage{theorem}
\usepackage{diagrams}
\theoremstyle{change}
\newtheorem{Thm}{Theorem}[section]
\newtheorem{Cor}[Thm]{Corollary}
\newtheorem{Prop}[Thm]{Proposition}
\newtheorem{Lem}[Thm]{Lemma}
{\theorembodyfont{\rmfamily}
\newtheorem{Num}[Thm]{}

}
\newcommand\wh{\widehat}
\newcommand\nb{\nobreakdash}
\def\={\overset{\text{def}}=}

\newcommand{\qed}{\ \hglue 0pt plus 1filll $\Box$}

\newcommand{\SKIP}[1]{}

\newcommand{\CC}{\mathbb{C}}   
\newcommand{\HH}{\mathbb{H}}

\newcommand{\RR}{\mathbb{R}}

\newcommand{\FF}{\mathbb{F}}
\newcommand{\OO}{\mathbb{O}}

\newcommand{\ZZ}{\mathbb{Z}}

\newcommand{\bra}[1]{\langle#1\rangle}
\newcommand{\into}{\hookrightarrow}

\newcommand{\ra}{\longrightarrow}

\newcommand{\BO}{\mathrm{BO}}
\newcommand{\isom}{\cong}

\newcommand{\Ahat}{\widehat{A}}

\DeclareMathAlphabet{\varcal}{U}{rsfs}{m}{it}
\begin{document}

\title{\bf A diffeomorphism classification
of manifolds which are like projective planes}
\author{Linus Kramer and Stephan Stolz%
\thanks{This paper was supported by a research grant 
in the {\sl Schwerpunktprogramm Globale Differentialgeometrie} by
the Deutsche Forschungsgemeinschaft; the second author was partially supported by NSF grant DMS 0104077.}}

\maketitle
\begin{abstract}
We give a complete
diffeomorphism classification of $1$-connected closed manifolds
$M$ with integral homology
$
H_*(M)\isom\ZZ\oplus\ZZ\oplus\ZZ
$,
provided that $\dim(M)\neq 4$.
\end{abstract}

The integral homology of an oriented closed manifold%
\footnote{All manifolds are assumed to be smooth.}
$M$ contains at least two copies of $\ZZ$ (in degree
$0$ resp.\ $\dim M$). If $M$ is simply connected and its homology has
minimal size (i.e., $H_*(M)\cong \ZZ\oplus \ZZ$), then $M$ is a homotopy
sphere (i.e., $M$ is homotopy equivalent to a sphere). It is well-known
from the proof of the (generalized) Poincar\'e conjecture
that any homotopy sphere is homeomorphic to the
standard sphere $S^n$ of dimension $n$. By contrast, the cardinality
of the set $\Theta_n$ of diffeomorphism classes of homotopy spheres
of dimension $n$ can be very large (but finite except possibly for $n=4$)
\cite{KM}. In fact, the connected sum of homotopy spheres gives $\Theta_n$ the structure of an abelian group which is closely related to the stable homotopy group $\pi_{n+k}(S^k)$, $k\gg n$ (currently known approximately in the range  $n<100$).

Somewhat surprisingly, it is easier to obtain an explicit diffeomorphism classification of $1$\nb-connected closed manifolds whose integral homology consists of {\em three} copies of $\ZZ$. Examples of such manifolds are the $1$\nb-connected projective planes (i.e., the projective planes over the complex numbers, the quaternions or the octonions). Eells and Kuiper pioneered the study of these `projective plane like' manifolds \cite{EK} and obtained many important and fundamental results. For example, they show that the integral cohomology ring of such a manifold $M$ is isomorphic to the cohomology ring of a projective plane, i.e., $H^*(M)\cong \ZZ[x]/(x^3)$. This in turn implies that the dimension of $M$ must be $2m$ with $m=2,4$ or $8$ (cf.\ \cite[\S 5]{EK}). 
We remark that a $1$\nb-connected closed manifold $M$ of dimension $n\ge 5$ with $H_*(M)\cong\ZZ\oplus\ZZ\oplus\ZZ$ admits a Morse function with three critical points, which is the assumption that Eells-Kuiper work with. Any $1$\nb-connected projective plane like manifold of dimension $4$ is homeomorphic to the complex projective plane by Freedman's homeomorphism classification of simply connected smooth $4$\nb-manifolds \cite{Fr}.

Eells and Kuiper prove that there are six (resp.\ sixty) homotopy types of projective plane like manifolds of dimension $2m$ for $m=4$ (resp.\ $m=8 $) \cite[\S 5]{EK}. They get close to obtaining a classification of these manifolds up to homeomorphism resp.\ diffeomorphism. One way to phrase their result is the following. If $M$ is a smooth manifold of this type, let $p^2_m(M)[M]\in\ZZ$ be the Pontryagin number obtained by evaluating the square of the Pontryagin class $p_m(M)\in H^{m}(M;\ZZ)$ (of the tangent bundle of $M$) on the fundamental class $[M]\in H_{2m}(M;\ZZ)$. Eells and Kuiper show that the Pontryagin number $p^2_m(M)[M]\in\ZZ$ determines the diffeomorphism type {\em up to connected sum with a homotopy sphere}; in other words, if $M'$ is another such manifold of the same dimension and the same Pontryagin number, then $M'$ is diffeomorphic to the connected sum $M\# \Sigma$ of $M$ with a $2m$\nb-dimensional homotopy sphere $\Sigma$ (see \ref{thm:classification} and \cite[\S 9]{EK}; we note that the Pontryagin number determines the Eells-Kuiper integer $h$ and vice versa via their formulas (2) resp.\ (5) in \S 9). 

A complete homeomorphism classification of topological manifolds
which look like projective planes was obtained by the
first author in \cite{K2}. The main result of this paper is the following.

\smallskip\noindent\textbf{Theorem A}
{\em Let $M$ be a smooth simply connected manifold of dimension $2m\ne 4$, with integral
homology $H_*(M)\cong\ZZ\oplus\ZZ\oplus\ZZ$. Then for any homotopy sphere
$\Sigma$ of dimension $2m$ the connected sum $M\#\Sigma$ is diffeomorphic to
$M$.}
\smallskip

In view of the results of Eells-Kuiper discussed above, this implies the following diffeomorphism classification of projective plane like manifolds.

\smallskip\noindent\textbf{Corollary B}
{\em
Let $M$ be a smooth simply connected $2m$\nb-manifold with integral
homology $H_*(M)\cong\ZZ\oplus\ZZ\oplus\ZZ$. Then the diffeomorphism type of $M$ is determined by the Pontryagin number $p^2_m(M)[M]\in\ZZ$.}
\smallskip

Results of Eells-Kuiper combined with a result of Wall \cite{Walln-1} allow a characterization of those integers which occur as the Pontryagin numbers of such manifolds. We will give a precise statement as Theorem \ref{thm:classification} in the next section; for now we remark that the above result provides us with an infinite family of manifolds $M$ which have a unique differentiable structure in the sense that any manifold homeomorphic to $M$ is in fact diffeomorphic to $M$ \cite{Kreck0}.
Another motivation for this paper came from the first author's attempt
to classify the underlying spaces of topological projective
planes in the sense of Salzmann \cite{CPP}. In 
\cite{K2} he obtained
a homeomorphism classification for the point sets of smooth
topological projective planes, showing that only the four
classical spaces $\FF\mathrm P^2$, $\FF=\RR,\CC,\HH,\OO$
appear. However, the diffeomorphism classification
remained open. Except for the case of $\CC\mathrm P^2$,
the results of the present paper settle this question.
Combininig our results with McKay's diffeomorphism classification of
$2$-dimensional smooth topological projective planes \cite{McKay},
we obtain the follwing result.

\smallskip\noindent\textbf{Corollary C}
{\em The point space of a smooth topological projective plane
(in the sense of \cite{CPP})
is diffeomorphic to its classical counterpart, i.e. to
$\RR\mathrm P^2$, $\CC\mathrm P^2$, $\HH\mathrm P^2$, or $\OO\mathrm P^2$.}

\medskip\textbf{Outline of the paper.}
In Section 1 we state in more detail the Eells-Kuiper results concerning the
diffeomorphism classification of projective plane like manifolds up to connected
sum with homotopy spheres. For the convenience of the reader, we also outline
the proofs. The other sections are devoted to proving our main Theorem A.
In Section 2 we use Kreck's modified surgery approach \cite{Kreck1} to show that for a closed simply connected manifold $M$ of dimension $2m\ne 4$ the connected sum $M\#\Sigma$ with a homotopy sphere $\Sigma$ is diffeomorphic to $M$ provided $\Sigma$ represents zero in a suitable bordism group $\Omega_{2m}^B$ (cf.\ Corollary \ref{cor}). 
The bordism groups $\Omega_*^B$ depend on a fibration $B\rTo BO$ which in turn depends on the manifold $M$. In Section 3 we determine the relevant fibration in the case that $M$ is a projective plane like $2m$\nb-manifold (cf.\ Proposition \ref{prop:normal_k-type}).  In Section 4 we prove that any homotopy sphere $\Sigma$ of dimension $2m=8,16$ represents zero in $\Omega_{2m}^B$ for $B$ as above, thus completing the proof of Theorem A.

\section{Classification up to connected sums with homotopy spheres}
As mentioned in the introduction, the diffeomorphism classification of 
projective plane like manifolds up to connected sum with homotopy spheres was obtained by Eells-Kuiper \cite{EK} (plus one result of Wall's \cite[Thm.\ 4, p.\ 178]{Walln-1}) or by specializing Wall's much more general classification of `almost closed' $(m-1)$\nb-connected $2m$ manifolds \cite{Walln-1} to this case. Still, we feel that it is worthwhile to outline in this section how this classification follows from the classification of $m$\nb-dimensional vector bundles over $S^m$ and the $h$\nb-cobordism theorem.

We recall that a smooth manifold $N$ is `almost closed' if it is a compact manifold whose boundary is a homotopy sphere. Such a manifold is obtained for example by removing an open $n$\nb-disk from a closed manifold $M$ of dimension $n$. The boundary $\partial N$ of an almost closed $n$\nb-manifold $N$ is homeomorphic to the standard sphere $S^{n-1}$, and we denote by $N(\alpha)=N\cup_\alpha D^n$ the closed topological manifold obtained by gluing $N$ and the disk $D^n$ along their common boundary via a homeomorphism $\alpha\colon \partial D^n\rTo \partial N$. We note that $N(\alpha)$ is again a smooth manifold if $\alpha$ is a {\em diffeomorphism}; moreover, if $\beta\colon \partial N\rTo \partial D^n$ is a second diffeomorphism, then $N(\beta)$ is diffeomorphic to the connected sum $N(\alpha)\#\Sigma$ of $N$ with the homotopy sphere 
$\Sigma=D^n\cup_{\alpha^{-1}\beta}D^n$ obtained by gluing two discs along their boundaries via the diffeomorphism $\alpha^{-1}\beta\colon \partial D^n\rTo \partial D^n$.

An almost closed manifold $N$ is called {\em projective plane like} if the integral homology $H_*(N(\alpha))$ (which is independent of the choice of the homeomorphism $\alpha$) is isomorphic to $\ZZ\oplus\ZZ\oplus\ZZ$. This implies that $N$ is a manifold of dimension $2m$ with $m=2,4,8$. 

\begin{Thm}[Eells-Kuiper] The diffeomorphism classes of simply connected almost closed projective plane like manifolds of dimension $2m$ for $m=4,8$ are in one-to-one correspondence with the non-negative integers. The manifold $N_t$ corresponding to $t\in\ZZ$ is the disk bundle of the vector bundle $\xi_t$ over $S^m$ with Euler class $e(\xi_t)= x$, and Pontryagin class $p_m(\xi_t)=2(1+2t)x$ (for $m=4$) resp.\ $p_m(\xi_t)=6(1+2t)x$ (for $m=8$), where $x$ is the generator of $H^m(S^m)$.
\end{Thm}

\begin{proof}
It is an easy homology calculation to show that the disk bundle $N_t=D(\xi_t)$ is an almost closed projective space like manifold (the condition $e(\xi_t)=x$ guarantees that the that the boundary $\partial D(\xi_t)$ is a homotopy sphere). We note that pulling $\xi_t$ back via a map $S^m\rTo S^m$ of degree $-1$ we obtain a bundle isomorphic to $\xi_{-t-1}$; it follows that the manifolds $N_t$ and $N_{-t-1}$   are diffeomorphic and hence it suffices to consider only
$t\geq 0$.

Conversely, if $N$ is any simply connected almost closed projective space like  manifold of dimension $2m$, consider the normal bundle $\xi$ of an embedding $S^m\hookrightarrow N$ which represents a generator for $H_m(N;\ZZ)\cong \ZZ$. 
Then $\xi$ is an $m$\nb-dimensional oriented vector bundle over $S^m$, whose disc bundle $D(\xi)$ can be identified with a tubular neighborhood of $S^m\subset N$. Up to isomorphism $\xi$ is determined by its Euler class $e(\xi)$ and its Pontryagin class $p_m(\xi)$. The assumption that $N$ is projective space like implies that the integral cohomology ring of $N/\partial N$ is isomorphic to $\ZZ[x]/(x^3)$, which in turn implies $e(\xi)=x$. By the classification of $m$\nb-dimensional vector bundles over $S^m$, this implies that $\xi$ is isomorphic to $\xi_t$ for some $t\in \ZZ$. Now removing the interior of the disc bundle $D(\xi)\subset N$ from $N$, we obtain a bordism $W$ between $\partial D(\xi)$ and $\partial N$. A homology calculation shows that this is in fact an $h$\nb-cobordism (i.e., the inclusion of either boundary component into $W$ is a homotopy equivalence). By Smale's $h$\nb-cobordism theorem, $W$ is diffeomorphic to $\partial D(\xi)\times [0,1]$; in particular, $N=D(\xi)\cup_{\partial D(\xi)} W$ is diffeomorphic to $D(\xi)=D(\xi_t)$, which proves the theorem.
\end{proof} 

The theorem above begs the question for which $t\in \ZZ$ is the boundary of $N_t$ {\em diffeomorphic} to the standard sphere $S^{2m-1}$. The answer is given by the next result:

\begin{Prop}[Eells-Kuiper, Wall] The boundary $\partial N_t$ is diffeomorphic to $S^{2m-1}$ if and only if  
$t\equiv 0,7,48,55\mod{56}$ (for $m=4$) resp.\ $t\equiv 0,127,16128,16255\mod{16256}$ (for $m=8$).
\end{Prop}

\begin{proof}
Choose a homeomorphism $\alpha\colon S^{2m-1}\rTo \partial N_t$ and consider the closed topological manifold $N_t(\alpha)=N_t\cup_\alpha D^{2m}$.
Its $\Ahat$\nb-genus $\Ahat(N_t(\alpha))$, a certain rational linear combination of the Pontryagin numbers $p_m^2(M)[M]$ and $p_{2m}(M)[M]$, turns out to be independent of $\alpha$, and can be expressed in terms of $t$ by the following formula \cite[\S 9, Thms.\ on p.\ 216 resp.\ p.\ 218]{EK}, \cite[\S 7.2]{K2}:
\[
\Ahat(N_t(\alpha))=
\begin{cases}
-\frac {t(t+1)}{7\cdot 8}& m=4\\
-\frac {t(t+1)}{127\cdot 128}& m=8
\end{cases}
\]

If $\partial N_t$ is diffeomorphic to $S^{2m-1}$, we may choose $\alpha$ to be a diffeomorphism, and then $N_t(\alpha)$ is a {\em smooth} manifold. This manifold can be equipped with a spin structure, since $H^i(M;\ZZ/2)=0$ for $i=1,2$ and hence the Stiefel-Whitney classes $w_i(M)\in H^i(M;\ZZ/2)$, $i=1,2$ (the potential obstructions against a spin structure) vanish. This implies that
$\Ahat(N_t(\alpha))$ is an {\em integer}, namely the index of the `Dirac operator' which can only be constructed for smooth spin manifolds. The formula above then implies that $t$ satisfies the congruence of the proposition. 

Conversely, according to a result of Wall \cite[Thm.\ 4, p.\ 178]{Walln-1}, the integrality of $\Ahat(N_t(\alpha))$ implies that $\partial N_t$ is 
diffeomorphic to the standard sphere.
\end{proof}
We note that the Pontryagin number $p_m^2(M)[M]$ of the projective plane like manifold $M=N_t(\alpha)$ is equal to $2^2(1+2t)^2$ (for $m=4$
--- unfortunately, the formula
stated in \cite{K2} p.~2 is off by a factor $2$;
the correct number given here
appears in \emph{loc.cit.} Thm.~7.1)
resp.\ $6^2(1+2t)^2$ (for $m=8$). Hence the theorem and the proposition above imply the following result.

\begin{Thm}\label{thm:classification}
 Let $M$ be a smooth projective plane like manifold of dimension $2m$, $m=4,8$. Then up to connected sum with a homotopy sphere, the diffeomorphism type of $M$ is determined by the Pontryagin number $p_m^2(M)[M]\in\ZZ$. Moreover, an integer $k$ is equal to the Pontryagin number $p^2_m(M)[M]$ of such a manifold if and only if $k$ is of the form $k=2^2(1+2t)^2$ with $t\equiv 0,7,48,55\mod{56}$ (for $m=4$) resp.\ $k=6^2(1+2t)^2$ with $t\equiv 0,127,16128,16255\mod{16256}$ (for $m=8$).
\end{Thm}

\section{Bordism groups and surgery}
In this section we briefly describe a main result of Kreck's `modified surgery theory' \cite{Kreck1} (Theorem  \ref{thm:Kreckthm} below). A direct consequence of this result (see Corollary \ref{cor}) is that the connected sum $M\#\Sigma$ of a closed simply connected manifold $M$ of dimension $2m\ne 4$ with a homotopy sphere $\Sigma$ is diffeomorphic to $M$ provided $\Sigma$ represents zero in a suitable bordism group $\Omega_{2m}^B$. 

We begin by defining the bordism groups $\Omega_n^B$. 

\begin{Num}
Fix a fibration $B\rTo BO$ over the classifying space $BO$ of the
stable orthogonal group. We recall that $BO$ is the union of the
classifying spaces $BO_k$ of the orthogonal groups and that $BO_k$
is the union of the Grassmann manifolds
$Gr_k(\RR^{n+k})$ of $k$-planes in
$\RR^{n+k}$ via natural inclusion
maps $Gr_k(\RR^{n+k})\subseteq Gr_{k}(\RR^{n+k+1})$.
Let $(N,\partial N)$ be a compact $n$-manifold, and let
$\iota:(N,\partial N)\rInto(\RR^{n+k}_+,\partial\RR^{n+k}_+)$
be a smooth embedding into euclidean half-space.
Recall that the normal Gauss map 
$\nu:N^n\rTo Gr_k(\RR^{n+k})$ assigns to any point
$x\in N$ its normal space in $\RR^{n+k}$.
A \emph{$B$-structure} on $N$ is an equivalence class of pairs
$(\iota,\bar\nu)$, where $\bar \nu$ is a map making the following 
diagram commutative 
\begin{diagram}[height=2em,width=4em]
&&&& B\\
&&&\ruTo(4,2)^{\bar\nu}
&\dTo\\
N & \rTo_{\nu} && Gr_k(\RR^{n+k})\subset BO_k\subset& BO.
\end{diagram}
The equivalence relation
is generated by simultaneous deformations of $\iota$ and $\bar\nu$,
and by the stabilization map $\RR^{n+k}\rInto\RR^{n+k+1}$.

A {\em $B$\nb-manifold} is a manifold equipped with a $B$\nb-structure; a {\em $B$\nb-bordism} between $B$\nb-manifolds $M_1\rTo^{\bar\nu_1}B$ and
$M_2\rTo^{\bar\nu_2}B$ is a bordism $W$ between $M_1$ and $M_2$ equipped with a $B$\nb-structure $\bar\nu\colon W\rTo B$ which restricts to $\bar\nu_1$ resp.\ $\bar\nu_2$ on the boundary $\partial W=M_1\cup M_2$. 
A $B$-structure $\beta\colon M\rTo B$ is called a \emph{normal $k$-smoothing}  if
$\bar\nu$ is a $(k+1)$-equivalence, i.e., if the induced homomorphism $\bar\nu_*\colon \pi_i(M)\rTo \pi_i(B)$ is an isomorphism for $i\le k$ and surjective for $i=k+1$. We remark that if there exists a $k$\nb-smoothing $\bar\nu\colon M\rTo BO$, the fibration $B\rTo BO$ is determined by the manifold $M$ up to fiber homotopy equivalence, if we assume that $\pi_iB\rTo \pi_iBO$ is an isomorphism for $i>k+1$ and injective for $i=k+1$. In this case, Kreck refers to $B\rTo BO$ as the {\em normal $k$\nb-type} of $M$.
\end{Num}

\begin{Thm}\label{thm:Kreckthm}
(Kreck \cite[Theorem B]{Kreck1} )
Let $M_1$, $M_2$ be closed manifolds of dimension $n\ge 5$ with the same Euler characteristic which are equipped with $B$\nb-structures that are normal $k$\nb-smoothings. For $k\geq [n/2]-1$ a $B$-bordism $W$ between $M_1$ and $M_2$  is bordant to an $s$-cobordism if and only if a certain obstruction $\theta(W)$ is elementary.
\end{Thm}

We recall that a bordism $W$ between $M_1$ and $M_2$ is an {\em $s$\nb-cobordism}  if the inclusions $M_1\rTo W$ and $M_2\rTo W$ are simple homotopy equivalences. The $s$-cobordism Theorem  implies 
  that then $M_1$ and $M_2$ are diffeomorphic (assuming that $\dim M_1=\dim M_2\ge 5$). 
\begin{Num}
The obstruction $\theta(W)$ is an element of an abelian monoid $\ell_{n+1}(\pi,w)$ which depends on the fundamental group $\pi=\pi_1(B)$ and the induced map $w\colon \pi_1(B)=\pi\rTo\pi_1(BO)=\ZZ/2$. Even if $\pi$ is the trivial group (this is the case we care about in this paper), the obstruction $\theta(W)$ is difficult to handle for $k=[n/2]-1$ (see \cite[\S 7]{Kreck1}). 
The situation greatly
simplifies for $k\ge [n/2]$: 
\begin{itemize}
\item The normal $k$\nb-smoothings induce isomorphisms $H_i(M_1)\cong H_i(B)\cong H_i(M_2)$ for $i\le  [n/2]$. By Poincar\'e duality, we also have isomorphisms $H_i(M_1)\cong H_i(M_2)$ for $ [n/2]+1\le i\le n$ and hence in particular the Euler characteristics of $M_1$ and $M_2$ agree.
\item By \cite{Kreck1} p.~734 the obstruction
$\theta(W)$
is contained in an abelian subgroup $L_{n+1}(\pi,w)$ of the
monoid $\ell_{n+1}(\pi,w)$. Moreover, this group 
projects to the Whitehead group $Wh(\pi)$, and the kernel
\[
L^s_{n+1}(\pi,w)=\mathrm{ker}\left(L_{n+1}(\pi,w)\rTo Wh(\pi)\right)
\]
is Wall's classical surgery group \cite{Wallsurgery}.
\end{itemize}
If $B$ is simply connected,
then $Wh(\pi)=0$, and   so $L^s_{n+1}(\pi,w)=L_{n+1}(\pi,w)$; moreover, these groups are zero if $n$  is even \cite{Wallsurgery}. As the zero-element in $L_{n+1}(\pi,w)$ is certainly
elementary in Kreck's sense, the obstruction
$\theta(W)$ is elementary in this case, and we conclude:
\end{Num}
\begin{Cor}
Let $M_1$, $M_2$ be closed simply connected $2m$\nb-dimensional manifolds which are equipped with $B$\nb-structures that are normal $m$\nb-smoothings, $m\ge 3$. If $M_1$ and $M_2$ represent the same element in the bordism group $\Omega_{2m}^B$, then $M_1$ is diffeomorphic to $M_2$.
\end{Cor}

\begin{Cor}\label{cor}
Let $\bar\nu\colon M\rTo B$ be a normal $m$\nb-smoothing of a simply connected $2m$\nb-ma\-ni\-fold, $m\ge 3$. Let $\Sigma$ be a homotopy sphere equipped with a $B$\nb-structure such that $[\Sigma]=0\in\Omega_{2m}^B$. Then $M\#\Sigma$ is diffeomorphic to $M$.
\end{Cor}

\begin{proof}
It is well-known that the connected sum $M\# N$ of two $B$\nb-manifolds admits a $B$\nb-structure such that it represents the same element in $\Omega_*^B$ as the disjoint union of $M$ and $N$; the desired $B$\nb-bordism $W$ is constructed by taking the disjoint union of $M\times [0,1]$ and $N\times [0,1]$ and attaching a $1$\nb-handle $D^n\times [0,1]$ to it connecting these two parts. The boundary of the resulting $n+1$\nb-manifold $W$ consists of the disjoint union of $M$, $N$ and $M\# N$; obstruction theory shows that the $B$\nb-structure can be extended over the $1$\nb-handle to give a $B$\nb-structure on $W$. 

We note that the $B$\nb-structure constructed on $M\#\Sigma$ in this way is again an $m$\nb-smoothing; hence the previous corollary implies that $M\#\Sigma$ is diffeomorphic to $M$. \qed\end{proof}

\section{The normal $m$-type of projective space like $2m$-manifolds}

In order to apply this result, we need to identify for a given projective plane like $2m$\nb-manifold $M$ a suitable fibration $B\rTo BO$ such that $M$ admits a normal $m$\nb-smoothing $\bar\nu \colon M\rTo B$ (i.e., $\bar\nu_*\colon \pi_iM\rTo \pi_iB$ is an isomorphism for $i\le M$ and surjective for $i=m+1$). To find $B$, we will need the following information about $M$.

\begin{Lem} \label{lem:induced} Let $M$ be a projective plane like $2m$\nb-manifold, $m=4,8$, such that the almost closed manifold $\mathring{M}$ obtained by removing an open disk from $M$ is diffeomorphic to $N_t$. Let $\nu\colon M\rTo BO_k\subset BO$ be the normal Gauss map induced by an embedding $M\subset \RR^{2m+k}$. Then the induced map $\nu_*\colon \pi_mM\cong\ZZ\ra \pi_mBO\cong\ZZ$ is multiplication by $\pm(2t+1)$.
\end{Lem}

\begin{proof}
Let $i\colon S^m\into D(\xi_t)=N_t\subset M$ be the inclusion of the zero-section.
The normal bundle of this embedding is $\xi_t$. The normal bundle of the embedding $N_t\subset M\into \RR^{n+k}$ is the pull back $\nu^*\gamma^k$ of the universal bundle $\gamma^k\rTo BO_k$ via the normal Gauss map $\nu\colon M\rTo BO_k$. This implies that the vector bundle
\[
TS^m\oplus \xi_t\oplus i^*\nu^*\gamma^k\cong
i^*TM\oplus  i^*\nu^*\gamma^k=i^*(TM\oplus \nu^*\gamma^k)
\]
is the restriction of the tangent bundle of $\RR^{n+k}$ to $S^m\subset M\subset \RR^{n+k}$ and hence trivial. Identifying stable vector bundles over $S^m$ with their classifying map $[S^m\rTo BO]\in \pi_mBO$, we conclude $i^*\nu^*\gamma=-\xi_t\in\pi_mBO$. Comparing the Pontryagin classes $p_m(\xi_t),p_m(\xi_1)\in H^m(S^m;\ZZ)$, we see that $\xi_t=(2t+1)\xi_1\in \pi_mBO$. Combining these facts, we have $i^*\nu^*\gamma=-(2t+1)\xi_1\in\pi_mBO$. 

Reinterpreting this equation, it tells us that the map $\nu_*\colon \pi_mM\rTo \pi_mBO$ maps the generator $[i\colon S^m\rTo M]\in\pi_mM$ to $-(2t+1)\xi_1\in \pi_mBO$, which implies the lemma, since $\xi_1$ is a generator of $\pi_mBO\cong\ZZ$.
\qed\end{proof}

\begin{Num}
We note that a projective plane like $2m$\nb-manifold $M$ is $(m-1)$\nb-connected, i.e., $\pi_iM=0$ for $i<m$. This implies by standard obstruction theory that the normal Gauss map $\nu\colon M\rTo BO_k\subset BO$ of an embedding $M\subset \RR^{2m+k}$ can be factored through the {\em $(m-1)$\nb-connected cover $q\colon BO\bra{m}\rTo BO$}, a fibration determined up to fiber homotopy equivalence by the requirement that $\pi_iBO\bra{m}=0$ for $i<m$ and that $q_*\colon\pi_iBO\bra{m}\rTo BO$ induces an  isomorphism for $i\ge m$ (we note that the $1$\nb-connected cover $X\bra{2}\rTo X$ of a space $X$ is just the universal covering of $X$). The lift $\bar\nu\colon M\rTo BO\bra{m}$ of $\nu$ constructed this way is {\em not} a normal $m$\nb-smoothing of $M$, since by the above lemma, the induced map $\pi_mM\rTo \pi_mBO\bra{m}=\pi_mBO$ is not an isomorphism unless $t=0$. In particular, $BO\bra{m}\rTo BO$ is not the normal $m$\nb-type of $M$ unless $t=0$.
\end{Num}

\begin{Num} Now we proceed to construct the fibration $B\rTo BO$ which will turn out to be the normal $m$\nb-type of $M$. Let $K(\ZZ,m)$ be the Eilenberg-MacLane space characterized up to homotopy equivalence by the requirement that the homotopy group $\pi_iK(\ZZ,m)$ is zero for $i\ne m$ and equal to $\ZZ$ for $i=m$. The long exact homotopy sequence of the path fibration 
\[
\Omega K(\ZZ,m)\ra P K(\ZZ,m)\ra K(\ZZ,m)
\]
together with the fact that the path space $PK(\ZZ,m)$ is contractible shows that the loop space $\Omega K(\ZZ,m)$ is the Eilenberg-MacLane space $K(\ZZ,m-1)$. For $m\equiv 0\mod{4}$, let us denote by $B_{d,m}\rTo BO\bra{m}$ the pull-back of the above path fibration via a map $\phi\colon BO\bra{m}\rTo K(\ZZ,m)$ such that the induced map $\pi_*\colon \pi_mBO\bra{m}=\ZZ\rTo \pi_mK(\ZZ,m)=\ZZ$ is multiplication by $d$ (this requirement determines $\phi$ up to homotopy). We note that the long exact homotopy sequence of this fibration shows that the induced map $\ZZ\cong\pi_mB_{d,m}\rTo \pi_mBO\bra{m}\cong\ZZ$
is multiplication by $\pm d$. 
\end{Num}

\begin{Prop}\label{prop:normal_k-type} Let $M$ be as in Lemma \ref{lem:induced}. Then the normal $m$\nb-type of $M$ is the composite fibration $B_{2t+1,m}\rTo BO\bra{m}\rTo  BO$. 
\end{Prop}

\begin{proof}
Let $\bar\nu'\colon M\rTo BO\bra{m}$ be the lift of the normal Gauss map $M\rTo BO$ associated to an embedding $M\into \RR^{n+k}$. Again obstruction theory shows that $\bar \nu'$ can be lifted to a map $\bar\nu\colon M\rTo B_{2t+1,m}$. Lemma \ref{lem:induced} implies that the induced map $\bar\nu_*\colon \pi_iM\rTo \pi_iB_{d,m}$ is an isomorphism for $i=m$. Moreover, it is surjective for $i=m+1$: for $m=4$ this is obvious, since $\pi_5BO=0$; for $m=8$, it follows from the fact that the Hopf map $\eta\colon S^9\rTo S^8$ induces a surjection $\pi_8BO=\ZZ\rTo\pi_9BO=\ZZ/2$.
\qed\end{proof}

Applying now Corollary \ref{cor} to projective plane like manifolds, we conclude:

\begin{Cor}
If $M$ is as in Lemma \ref{lem:induced}, and $\Sigma$ is a homotopy sphere of dimension $2m$ with $[\Sigma]=0\in \Omega_{2m}^{B_{2t+1,m}}$, then $M\#\Sigma$ is diffeomorphic to $M$.
\end{Cor}

\section{The bordism class of homotopy spheres}

In view of the last corollary our main result follows from the following statement whose proof 
is the goal of this section.

\begin{Prop}
Let $\Sigma$ be a homotopy sphere of dimension $2m=8,16$. Then 
$[\Sigma]=0\in \Omega_{2m}^{B_{2t+1,m}}$ for any $t$.
\end{Prop}

To prove this result we note that the map $B_{2t+1,m}\rTo BO\langle m\rangle $ is a map of fiber bundles over $BO$ and hence it induces a homomorphism of bordism groups
\begin{equation}\label{eq:bordismmap}
\Omega_{*}^{B_{2t+1,m}}\ra \Omega_{*}^{BO\langle m\rangle}.
\end{equation}
The groups $\Omega_{2m}^{\BO\bra{m}}$ are known for $m=4,8$,
see Milnor \cite{Mil} and Giambalvo \cite{Gam}:
\[
\Omega_8^{BO\bra4}=\Omega_8^{BSpin}\cong\ZZ\oplus\ZZ
\quad\text{ and }\quad
\Omega_{16}^{BO\bra{8}}\cong\ZZ\oplus\ZZ.
\]
Since $\Theta_{2m}\cong \ZZ/2$ for $m=4,8$, it follows that for any homotopy $2m$\nb-sphere $\Sigma$ the connected sum $\Sigma\#\Sigma$ is diffeomorphic to $S^{2m}$. In particular, $\Sigma$ represents an element of order at most
$2$ in $\Omega_{*}^{B_{2t+1,m}}$. Hence the next result implies the proposition above.

\begin{Lem} The homomorphism \eqref{eq:bordismmap} is a $2$\nb-local isomorphism (i.e., its kernel and and cokernel belong to the class of torsion groups  without elements of order $2$).
\end{Lem}

Before proving this lemma we recall some relevant facts.

\begin{Num}{\bf The Pontryagin-Thom construction.}
Let $f_k\colon B_k\rTo BO_k$ be the restriction of the fibration $f\colon B\rTo BO$ to $BO_k\subset BO$. Let $\gamma_k\rTo BO_k$ be the universal $k$\nb-dimensional vector bundle, let $\wh \gamma^k\rTo B_k$ be its pull-back via $f_k$, and let $T(\wh \gamma_k)$ be the Thom space of $\wh \gamma^k$ (the quotient space of its total space obtained by collapsing all vectors of length $\ge 1$ to a point). Then the Pontryagin-Thom construction (see \cite[Thm\., p. 18]{Stong} produces an isomorphism
\begin{equation}\label{eq:PTisom}
\Omega_n^B\cong \lim_{k\to\infty}\pi_{n+k}T(\wh \gamma^k).
\end{equation}
\end{Num}

\begin{Num}{\bf Thom spectra.}
It is usual and convenient to express the right hand side of the Pontryagin-Thom isomorphism \eqref{eq:PTisom} in terms of Thom spectra. We recall that a {\em spectrum} is a sequence $E_k$ of pointed spaces together with pointed maps $\Sigma E_k\rTo E_{k+1}$ from the suspension of $E_k$ to $E_{k+1}$. For example, if $B\rTo BO$ is a fibration, there is an associated {\em Thom spectrum} $M\!B$, whose $k$\nb-th space is the Thom space $T(\wh \gamma^k)$.

Many constructions with spaces can be generalized to spectra; e.g., the homotopy (resp.\ homology) groups of a spectrum $E=\{E_k\}$ are defined as 
\[
\pi_nE\=\lim_{k\to\infty}\pi_{n+k}E_k
\qquad
H_n(E)\=\lim_{k\to\infty}H_{n+k}(E_k).
\]
With these definitions, the Pontryagin-Thom isomorphism takes the pleasant form
\begin{equation}\label{eq:PTspectra}
\Omega_n^B\cong \pi_n M\!B.
\end{equation}
Assuming that $B$ is $1$\nb-connected, the vector bundles $\wh\gamma^k\rTo B_k$ are all oriented and hence we have Thom-isomorphisms $H_i(B_k;\ZZ)\cong H_{i+k}(T\wh\gamma^k;\ZZ)$. It turns out that these isomorphisms are all compatible and so passing to the $k\to\infty$ limit,  one obtains a Thom-isomorphism
\begin{equation}\label{eq:Thomisom}
H_i(B;\ZZ)\cong H_i(M\!B;\ZZ)
\end{equation}
\end{Num}

\par\medskip\rm\emph{Proof of lemma.}
By construction, the induced homomorphism
\[
\pi_i(B_{2t+1,4k})\ra \pi_i (BO\bra{4k})
\]
is an isomorphism for $i\ne 4k$; for $i=4k$, it is injective with cokernel isomorphic to $\ZZ/(2t+1)$. In particular, for all $i$ it is a $2$\nb-local isomorphism. Then the generalized Whitehead Theorem (\cite[Thm.\ 22 of Chap.\ 9, \S 6]{Span}) implies that 
\begin{equation}\label{eq:p}
p_*\colon H_i(B_{2t+1,4k};\ZZ)\ra H_i (BO\bra{4k};\ZZ)
\end{equation}
is also a $2$\nb-local isomorphism. Now we consider the map of Thom spectra 
\[
Mp\colon M\!B_{2t+1,4k}\rTo M\!BO\bra{4k}
\]
induced by $p$. Via the Thom isomorphism \eqref{eq:Thomisom} the induced map in {\em homology} may be identified with the homomorphism \eqref{eq:p}, while the induced map on {\em homotopy groups} via the Pontryagin-Thom isomorphism \eqref{eq:PTspectra} corresponds to the homomorphism \eqref{eq:bordismmap} of bordism groups. Again by the generalized Whitehead Theorem, the latter is a $2$\nb-local isomorphism since the former is.
\qed

\bibliographystyle{plain}

\begin{thebibliography}{WWW}

\bibitem{AdamsStable}
J. F. Adams,
{\it Stable homotopy and generalised homology},
Reprint of the 1974 original, Univ. Chicago Press,
Chicago, IL, 1995.
MR1324104 (96a:55002)
Zbl 0309.55016

\bibitem{AA}
J. F. Adams\ and\ M. F. Atiyah,
$K$-theory and the Hopf invariant,
Quart. J. Math. Oxford Ser. (2) {\bf 17} (1966), 31--38.
MR0198460 (33 \#6618)
Zbl 0136.43903

\bibitem{Bott}
R. Bott,
The stable homotopy of the classical groups,
Ann. of Math. (2) {\bf 70} (1959), 313--337.
MR0110104 (22 \#987)
Zbl 0129.15601

\bibitem{EK} 
J. Eells, Jr.\ and\ N. H. Kuiper,
Manifolds which are like projective planes,
Inst. Hautes \'Etudes Sci. Publ. Math. No. 14 (1962), 5--46.
MR0145544 (26 \#3075)
Zbl 0109.15701

\bibitem{Fr} 
M. H. Freedman,
The topology of four-dimensional manifolds,
J. Differential Geom. {\bf 17} (1982),
no.~3, 357--453.
MR0679066 (84b:57006)
Zbl 0528.57011

\bibitem{Gam}
V. Giambalvo,
On $\bra8$-cobordism, Illinois J. Math.
{\bf 15} (1971), 533--541.
MR0287553 (44 \#4757)
Zbl 0221.57019

V. Giambalvo, Correction to my paper:
``On $\bra8$-cobordism'',      
Illinois J. Math. {\bf 16} (1972), 704.
MR0309132 (46 \#8243)
Zbl 0238.57020

\bibitem{KM}
M. A. Kervaire\ and\ J. W. Milnor,
Groups of homotopy spheres. I,
Ann. of Math. (2) {\bf 77} (1963), 504--537.
MR0148075 (26 \#5584)
Zbl 0115.40505

\bibitem{K1}
L. Kramer,
The topology of smooth projective planes,
Arch. Math. (Basel) {\bf 63} (1994),
no.~1, 85--91.
MR1277915 (95k:51021)
Zbl 0831.51011

\bibitem{K2}
L. Kramer,
Projective planes and their look-alikes,
J. Differential Geom. {\bf 64} (2003),
no.~1, 1--55.
MR2015043 (2004g:57041)
Zbl 1068.57019

\bibitem{Kreck0}
M. Kreck,
Manifolds with unique differentiable structure,
Topology {\bf 23} (1984), no.~2, 219--232
MR0744852 (85j:57051)
Zbl 0547.57025

\bibitem{Kreck1}
M. Kreck,
Surgery and duality,
Ann. of Math. (2) {\bf 149} (1999), no.~3, 707--754.
MR1709301 (2001a:57051)
Zbl 0935.57039

\bibitem{Kreck2}
M. Kreck,
A guide to the classification of manifolds,
in {\it Surveys on surgery theory, Vol. 1},
121--134, Ann. of Math. Stud., 145,
Princeton Univ. Press, Princeton, NJ.
MR1747533 (2001b:57061)
Zbl 0947.57001

\bibitem{MM} 
I. Madsen\ and\ R. J. Milgram,
{\it The classifying spaces for surgery and cobordism of manifolds},
Ann. of Math. Stud., 92,
Princeton Univ. Press, Princeton, N.J., 1979.
MR0548575 (81b:57014)
Zbl 0446.57002

\bibitem{McKay}
B. McKay,
Smooth projective planes,
Geom. Ded  {\bf 116} (2005), 157--202.
MR2195446
Zbl 05021438

\bibitem{Mil}
J. Milnor,
Spin structures on manifolds,
Enseignement Math. (2) {\bf 9} (1963), 198--203.
MR0157388 (28 \#622)
Zbl 0116.40403

\bibitem{MS} 
J. W. Milnor\ and\ J. D. Stasheff,
{\it Characteristic classes},
Ann. of Math. Stud., 76,
Princeton Univ. Press, Princeton, N. J., 1974.
MR0440554 (55 \#13428)
Zbl 0298.57008

\bibitem{Rud}
Y. B. Rudyak,
{\it On Thom spectra, orientability, and cobordism},
Springer, Berlin, 1998. 
MR1627486 (99f:55001)
Zbl 0906.55001

\bibitem{CPP} 
H. Salzmann\ et al.,
{\it Compact projective planes},
de Gruyter, Berlin, 1995.
MR1384300 (97b:51009)
Zbl 0851.51003

\bibitem{Span}
E. H. Spanier,
{\it Algebraic topology},
Corrected reprint, Springer, New York, 1981.
MR0666554 (83i:55001)
Zbl 0477.55001

\bibitem{Stolz}
S. Stolz,
{\it Hochzusammenh\"angende Mannigfaltigkeiten und ihre R\"ander},
Springer, Berlin, 1985.
MR0871476 (88f:57061)  
Zbl 0561.57021

\bibitem{Stong}
R. E. Stong,
{\it Notes on cobordism theory}, 
Mathematical notes 
Princeton University Press, Princeton, N.J.;
University of Tokyo Press, Tokyo 1968
MR 0248858 (40 \#2108)
Zbl 0181.26604

\bibitem{Walln-1}
C. T. C. Wall,
Classification of $(n-1)$-connected $2n$-manifolds.
Ann. of Math. (2)  75  1962 163--189.
MR0145540 (26 \#3071)
Zbl 0218.57022
 
\bibitem{Wallsurgery} 
C. T. C. Wall,
{\it Surgery on compact manifolds},
Second edition, Amer. Math. Soc., Providence, RI, 1999.
MR1687388 (2000a:57089)
Zbl 0935.57003

\end{thebibliography}

\bigskip
\raggedright
Linus Kramer\\
Mathematisches Institut, 
Universit\"at M\"unster,
Einsteinstr. 62,
48149 M\"unster,
Germany\\
\makeatletter
e-mail: {\tt linus.kramer{@}math.uni-muenster.de}\\\smallskip
Stephan Stolz\\
Department of Mathematics, University of Notre Dame, Notre Dame, IN 46556,
USA\\
e-mail: {\tt stolz.1{@}math.nd.edu}

\end{document}